\documentclass{amsart}

\usepackage{amsmath}
\usepackage{amssymb}

\newcommand{\dd}{\textrm{d}}
\newcommand{\CC}{\mathbb{C}}

\title{An integral representation for the Lambert $W$ function}
\author{Istv\'an Mez\H{o}}
\address{Department of Mathematics\\Nanjing University of Information Science and Technology\\Nanjing, P. R. China}
\email{istvanmezo81@gmail.com}

\begin{document}

\begin{abstract}
Based on a Problem and its solution published on the pages of SIAM Review, we give an interesting integral representation for the Lambert $W$ function in this short note. In particular, our result yields a new integral representation for the $\Omega=W(1)$ constant as well.\end{abstract}

\maketitle

\section{Statement of the result}
Thirty five years ago A. Nuttall \cite{Nuttall} posed a problem on the pages of SIAM Review, and C. Bouwkamp \cite{Bouwkamp} find the solution a year later:
\begin{equation}
\int_{0}^{\pi}\left[\frac{\sin t}{t} \exp (t \cot t)\right]^{\nu} d t=\frac{\pi \nu^{\nu}}{\Gamma(1+\nu)} \quad \text { for } \nu \geq 0.\label{NB_int}
\end{equation}
Here $\nu$ is real number, and $\Gamma$ is the Euler Gamma function.

Based on this result, we give a new integral evaluation for the Lambert $W$ function which is becoming more-and-more important in applied mathematics \cite{W}. The $W(x)$ function is defined as the solution(s) of the transcendental equation
\[we^w=x.\]
This is a many valued function, for which two branches, denoted by $W_{-1}$ and $W_0$ take real values. $W_0$ is the principal branch, and it is represented around $x=0$ by the Taylor series
\begin{equation}
W(x)=\sum_{n=1}^\infty\frac{(-n)^{n-1}}{n!}x^n\quad\left(|x|<\frac1e\right).\label{id:basicprop_WTaylor}
\end{equation}

Our result is that the following representation is valid for the principal branch $W_0(x)$:
\begin{equation}
W_0(x)=\frac{1}{\pi}\int_0^\pi\log\left(1+x\frac{\sin t}{t}e^{t\cot t}\right)\dd t.\label{id:basicprop_2ndintrep}
\end{equation}

The value $\Omega=W_0(1)$ has a special role \cite{Finch}. Setting $x=1$ in the above formula, we get
\begin{equation}
\Omega=\frac{1}{\pi}\int_0^\pi\log\left(1+\frac{\sin t}{t}e^{t\cot t}\right)\dd t.\label{Omega}
\end{equation}

\section{Proof}

We start with the Nuttall\,--\,Bouwkamp integral. We only need the fact that $\Gamma(\nu+1)=\nu!$ for non-negative a integer $\nu$. Divide both sides of \eqref{NB_int} by $\nu$, multiply with $(-1)^{\nu-1}x^\nu$, and sum for all $\nu\ge1$ integers. On the left-hand-side we get
\[\int_0^\pi\log\left(1+x\frac{\sin t}{t}e^{t\cot t}\right)\dd t.\]
On the right-hand-side we have the following sum:
\[\pi\sum_{\nu=1}^\infty\frac{(-\nu)^{\nu-1}}{\nu!}.\]
The sum is nothing else but the \eqref{id:basicprop_WTaylor} Taylor series of the principal branch of the Lambert $W$ function. Our result is therefore proven.

Although we used the Taylor series to connect the integral to $W_0$, numerical calculations show that \eqref{id:basicprop_2ndintrep} is valid on a large subset of $\CC$. In particular, although the point $x=1$ does not belong to the convergence domain of the Taylor series, \eqref{Omega} is still valid. The exact set on which \eqref{id:basicprop_2ndintrep} is valid is not known and can be an interesting problem in itself.

\section*{Acknowledgement}

After the previous version of the manuscript was prepared, G. Nemes pointed out, that the representation \eqref{id:basicprop_2ndintrep} is actually not new, it was known by G. A. Kalugin et al. \cite{KJC} since 2011. We would like to thank G. Nemes for this information.

\end{document}